\documentclass[12pt]{article}
\usepackage{amssymb}
\newtheorem{thm}{Theorem}[section]
\newtheorem{dfn}{Definition}[section]
\newtheorem{prop}{Proposition}[section]
\newtheorem{lem}{Lemma}[section]
\newtheorem{rem}{Remark}[section]

\newtheorem{prob}{Problem}[section]

\newtheorem{ex}{Example}[section]

\def\b{{\mathfrak{B}}}
\def\c{{\mathfrak{C}}}
\def\s{{\mathfrak{S}}}
\def\D{{\mathfrak{D}}}
\def\g{{\mathfrak{G}}}
\def\H{{\mathfrak{H}}}
\def\HH{{\mathbb{H}}}

\def\Z{{\mathbb{Z}}}

\def\Q{\mathbb{Q}}

\def\Z{\mathbb{Z}}
\def\qed{{\blacksquare}}
\begin{document}

{}

\vspace{2cm}

\begin{center}
{\Large \bf On Double Schubert and Grothendieck polynomials for Classical 
Groups }
\end{center}

\vspace{0.5cm}

\begin{center}
{\bf  A.N.Kirillov}
\end{center}

\begin{center}
Research Institute of Mathematical Sciences, RIMS, \\
Kyoto University, Sakyo-ku, 606-8502, Japan \\
{\it URL: ~~~~http://www.kurims.kyoto-u.ac.jp/\kern-.05cm$\tilde{\quad}$\kern-.17cm kirillov }
$$and$$
{The Kavli Institute for the Physics and Mathematics of the Universe 
( IPMU ),\\ 5-1-5 Kashiwanoha,  Kashiwa, 277-8583, Japan} \\
\end{center}

\vspace{1cm}

\vspace{1cm}
\begin{abstract}
We give an algebra-combinatorial constructions of (noncommutative) 
generating functions of double Schubert and double $\beta$-Grothendieck 
polynomials corresponding to the full flag varieties associated to the Lie 
groups of classical types $A,B, C$ and $D$. Our approach is based on the 
decomposition of certain `` transfer matrices `` corresponding to the exponential solution to the quantum Yang--Baxter equations associated with either NiCoxeter or IdCoxeter salgebras of classical types. 

The ``triple''~$\beta$-Grothendieck polynomials ${\mathfrak{G}}_{w}^{W}(X,Y,Z)$ we have introduced, satisfy, among other things, the coherency and 
(generalized) vanishing conditions. Their generating function has a nice 
factorization in the  algebra $Id_{\beta}Coxeter(W)$, and as a 
consequence, 
the polynomials ${\mathfrak{G}}_{w}^{W}(X,Y,Z)$ admit a combinatorial description in terms of  $W$-type pipe dreams.

\end{abstract}
\section{Introduction}
Let $G$ be a Lie group of one of the classical types $A_{n-1},B_{n},C_{n},
D_{n}.$~Let $B$ be a Borel subgroup, $B^{-}$ be the opposite Borel subgroup, 
$T=B \cap B^{-}$ be the  maximal torus and $W:=W(G)$ be the Weyl group.~ 
Let ${\cal X}=G/B$ be the flag variety of a classical type. The description of 
the equivariant cohomology ring $H^{*}_{T}(G/B,\Z)$ of the flag variety $G/B$ 
is well-known, and can be presented in the form  
$$H^{*}_{T}(G/B,\Z)=\Z[x_1,\ldots,x_n, y_1, \ldots,y_n] /J_n,$$ 
where $J_n$ is the ideal generated by

$(type ~{\bf A_{n-1}})$~~~$e_i(x_1,\ldots,x_n)-e_i(y_1,\ldots,y_n),$~
$ 1\le i \le n.$ 

$(type ~{\bf C_n}~ or ~{\bf B_n})$~~~$e_i(x_1^2,\ldots,x_n^2)-e_i(y_1^2,\ldots,
y_n^2),$~$ 1\le i \le n.$ 

$(type {\bf D_n})$~~~$e_i(x_1^2,\ldots,x_n^2)-e_i(y_1^2,\ldots,y_n^2),$~~
$ 1\le i \le n-1,$ 

~~~~~~~~~~~~~~~~~~~and ~$e_n(x_1,\ldots,x_n)-e_n(y_1,\ldots,y_n).$

There are distinguish elements $[{\cal X}_w]_{T}$ in the cohomology 
ring $H^{*}_{T}(G/B,\Z),$ namely, the Poincare dual classes of the homology 
classes corresponding to the Schubert subvarieties 
${\cal X}_{w}={\overline {B~w~B/B}} \subset {\cal X}.$ ~In the cohomology 
ring $H^{*}_{T}(G/B,\Z)$ one can write 
$[{\cal X}_w]:= {\cal X}_{w}(X_n,Y_n)$ for a certain (homogeneous) polynomial 
${\cal X}_{w}(X_n,Y_n)$ of degree $l(w)$ in each set of variables $X_n$ and 
$Y_n.$  The sets of variables $X_n:=(x_1, \ldots, x_n)$ and 
$Y_n:=(y_1, \ldots,y_n)$ are known as the {\it Borel generators} of the 
equivariant cohomology ring $H^{*}_{T}(G/B,\Z);$ the variables $Y_n$ 
correspond to generators of the 
equivariant cohomology ring of a point, $H_{T}^{*}(pt)=\Z[Y_n]$, and the set 
of variables $X_n$ comes from the Chern classes of some linear vector bundles 
over the flag variety in question. By definition, the equivariant Schubert 
polynomials, or double Schubert polynomials ${\cal X}_{w}(X_n,Y_n),$ are 
polynomials which express the equivariant classes $[{\cal X}_{w}]_{T}$ in 
terms of Borel generators.

Polynomials ${\cal X}_{w}(X_n,Y_n),$ ~$ w \in W,$ 
are defined only modulo the ideal $J_n$. It is known that for any 
finite dimensional, semisimple Lie group of rank $n$ the set of 
polynomials ${\cal X}_{w}(X_n,Y_n),$ ~$ w \in W,$ possess and characterized 
by the following properties ( modulo the ideal $J_n$ of relations in the 
cohomology ring $H^{*}_{T}(G/B,\Z)$)

$({\bf A})$~~~Polynomials ${\cal X}_{w}(X_n,Y_n),$~$w \in W,$ form a 
$\Z[Y]$-linear basis of the cohomology ring $H^{*}_{T}(G/B,\Z);$

$({\bf B})$~~~Coherency conditions)~~For any simple root $\alpha,$
$$\partial_{\alpha}^{(x)}~{\cal X}_{w}(X_n,Y_n)= \left\{\matrix{
{\cal X}_{ws_{\alpha}}(X_n,Y_n),~~l(ws_{\alpha})=l(w)+1, \cr
0~~~~~otherwise; \cr} \right.$$
$$\partial_{\alpha}^{(y)}~{\cal X}_{w}(X_n,Y_n)= \left\{\matrix{
{\cal X}_{s_{\alpha}w}(X_n,Y_n),~~l(s_{\alpha}w)=l(w)+1, \cr
0~~~~~otherwise; \cr} \right.$$

$({\bf C})$~~~(Vanishing conditions)~~Let $w,v \in W,$ then
$${\cal X}_{w}(X_n,-v(X_n)) =0,~~unless~~v \le w,$$
where the symbol $\le$ denotes the Bruhat order on the group $W;$ 

$({\bf D})$~~~(Normalization condition)~~${\cal X}_{id}(X_n,Y_n) =1,$ where 
$id \in W$ is the identity element.

Recall that $s_{\alpha}$ stands for the reflection 
corresponding to a simple root ${\alpha};$ ~~$\partial _{\alpha}^{(x)}=
{1-s_{\alpha} \over \alpha}$ denotes the corresponding Demazure operator 
acting on the variables $X_n;$~$v \in W$ acts on $X_n$ via the reflection 
representation.

In the case of flag varieties corresponding to the classical groups, 
polynomials   ${\cal X}_{w}(X_n,Y_n)$ possess also the so-called 
{\it stability property:}

$({\bf E})$~~Let $G_n$ be a Lie group of one of the classical series, and 
$\iota : G_n \hookrightarrow G_{n+1}$ be the canonical inclusion corresponding to the Dynkin diagram's embedding. If $w \in G_n,$ then
$$ {\cal X}_{w}(X_n,Y_n)=
{\cal X}_{\iota(w)}(X_{n+1},Y_{n+1})_{x_{n+1}=0=y_{n+1}}.$$
In the case of type $A_{n-1}$ flag varieties,  A. Lascoux and 
M.-P. Sch\"{u}utzenberger have constructed a family of polynomials 
$\{ \s_{w}(X_n,Y_n) \in \Z_{ \ge 0}[X_n,Y_n], w \in {\mathbb S}_n \},$ 
called double Schubert polynomials, that satisfies the all properties 
$({\bf A})$--$({\bf E})$ listed above, see e.g. \cite{Mac},\cite{Man} for 
detail account. It happened that the Lascoux--Sch\"{u}tzenberger double 
Schubert polynomials have non-negative integer coefficients and possess many 
nice combinatorial and algebraic properties. One of the basic properties of 
the double Schubert polynomials is the Cauchy type identity, that connects the 
simple Schubert polynomials $\s_{w}(X_n):=\s_{w}(X_n,0)$ with the double ones. 
Namely,
\begin{equation}
\s_{w}(X_n,Y_n)=\sum_{u,v}~\s_{u}(X_n)~\s_{v}(Y_n),
\end{equation}
where the sum runs over all $u,v \in {\mathbb S}_n$ such that $w=v^{-1}~u$
and $l(w)=l(u)+l(v).$ 

Recall that $\s_{w}(X_n) \in \Z_{\ge 0}[X_n]$ denotes 
the Lascoux-- Sch\"{u}tzenberger 
Schubert polynomial corresponding to a permutation $w \in {\mathbb S}_n.$ The 
set of Schubert polynomials $\s_{w}(X_n), w \in {\mathbb S}_n$ satisfy the 
Stability and Coherency Conditions (without passing to the quotient modulo 
the ideal $J_n$ !)  and their images in the cohomology ring 
$H^{*}({\cal F}l_n,\Z)$ form a basis. Conversely, if $W$ is a Weyl group of 
a classical type, and one has a family of
polynomials $\phi_{w}(X_n) \in \Z[X_n], w\in W,$ which satisfies the 
conditions $({\bf A})--({\bf E})$ (except that ${(\bf C})$), then the family 
of double polynomials  
\begin{equation} 
\Phi_{w}(X_n,Y_n)=\sum_{u,v}~\phi_{u}(X_n)~\phi_{v}(Y_n),
\end{equation}
summed over all $u,v \in W $ such that $w=v^{-1}~u$ and $l(w)=
l(u)+l(v),$ {\it also satisfies} the conditions $({\bf A})$--$({\bf E}),$ 
except, probably, the Vanishing Conditions. This brings up the natural 
question: 

Consider any Weyl group $W$ of a classical type, does there exist a 
family of polynomials $\phi_{w}(X_n) \in \Z[X_n], w\in W$ such that the family 
of double polynomials $\Phi_{w}(X_n,Y_n),~w \in W$ defined by $(2)$ satisfies 
the Vanishing Conditions $({\bf C})$ ?.

The main goal of the present paper is to show that the answer on this 
question is {\it Yes}. Namely, the Schubert polynomials of the second kind 
introduced originally in \cite{BH} for all Weyl groups of classical types, 
generate the set of double polynomials, called $(B,C,D)$-double Schubert 
polynomials of the second kind,  that satisfy the Vanishing conditions 
(without passing to the quotient). As it was observed in \cite{FK3}, the 
Schubert polynomials of the first kind introduced in \cite{FK3}, are closely 
related with the Schubert polynomials for classical groups 
introduced by S. Billey and M. Haiman \cite{BH}. Thus, the main result of the 
present work can be considered as a generalization of the construction of 
Schubert polynomials for classical groups given in \cite{BH} and \cite{FK3},  
to the case of double Schubert polynomials. We also give a construction of 
double Grothendieck polynomials for classical groups. Our main result in the 
case of Schubert polynomials can be stated as follows.

Let $W$ be a Weyl group of one the types $B,C$ or $D.$ Consider the 
Nil--Coxeter algebra $Nil(W),$ see Section~2.3, and the corresponding elements 
$B^{W}(x) \in \Z[x][Nil(W)].$ ~Finally, for any $w \in W$ define {\it double 
Schubert polynomial} $\s^{W}_{w}(X_n,Y_n)$ via the decomposition
$$(\s^{A_{n-1}}(-Y_n))^{-1}~\sqrt{B^{W}(Y_n)B^{W}(-X_n)}~\s^{A_{n-1}}(X_n)=
\sum_{w \in W}~\s^{W}_{w}(X_n,Y_n)~u_{w},$$
where $\s^{A_{n-1}}(X_n)$ denotes the Schubert expression of type $A_{n-1},$
introduced in \cite{FS}, ~and $B^{W}(X_n)=B^{W}(x_1) \cdots B^{W}(x_n).$
\begin{thm}~~~The set of  polynomials 
$\s^{W}_{w}(X_n,Y_n) \in \Z[X_n,Y_n],~w \in W$ satisfies the all conditions 
$({\bf A})$--$({\bf E}),$ and can be taken as a system of representatives for 
the equivariant Schubert classes $[{\cal X}_{w}]_{T} \in H^{*}_{T}(G/B,\Z).$
\end{thm}
We define double Grothendieck polynomials $\g_{w}^{W}(X_n,Y_n)$ corresponding 
to a Weyl group $W$ of classical type 
\footnote{in the case of Weyl groups of 
type $C$ one has to use the function $\phi_{2\beta}(-X_n)$ instead of that 
$\phi_{\beta}(-X_n).$ }
via the decomposition of the expression
$$(\s^{A_{n-1}}(\phi_{\beta}(-Y_n))^{-1}~
\sqrt{{\cal B}^{W}(Y_n){\cal B}^{W}(X_n)}
~\s^{A_{n-1}}(\phi_{\beta}(X_n))=
\sum_{w \in W}~\g^{W}_{w}(X_n,Y_n)~u_{w}$$
in the Id--Coxeter algebra $Id_{\beta}(W).$ 
~Hereinafter,  $\phi_{\beta}(x):= x/(1-{\beta \over 2}~x)$ ~and~
$\phi_{\beta}(-X_n)=(\phi_{\beta}(-x_1), \ldots, \phi_{\beta}
(-x_n)).$
\begin{thm}~~The double Grothendieck polynomials $\g_{w}^{W}(X_n,Y_n)$ 
corresponding to a Lie group of classical type, satisfy the all conditions 
$({\bf A})$--$({\bf E}),$ if one replaces the divided difference operators in 
the Coherency conditions $({\bf B}),$ on the ${\underline {isobaric}}$ divided 
difference operators. 
\end{thm}

To study combinatorial properties of the double Schubert polynomials 
$\s^{W}_{w}(X_n,Y_n)$ of the second kind, as well as to reveal their 
connections with the polynomials introduced and studied by S. Billey and 
M. Haiman \cite{BH}, it is convenient to introduce the set of polynomials 
$\s^{W}_{w}(X_n,Y_n, Z_m)$ depending on three set of variables via the 
decomposition in the Nil--Coxeter algebra $Nil(W)$ of the following expression:$$(\s^{A_{n-1}}(-Y_n))^{-1}~{B^{W}(Z_m)}~\s^{A_{n-1}}(-X_n)=
\sum_{w \in W}~\s^{W}_{w}(X_n,Y_n,Z_m)~u_{w}.$$ 
These polynomials are common generalization of both the Stanley symmetric 
polynomials $F_{w}^{W}(Z_{m})$ of type $W,$ coming from the decomposition 
${B^{W}(Z_m)}=\sum_{w \in W}F_{w}(Z_m)~u_w,$ and the double Schubert 
polynomials of the first kind introduced in \cite{FK3} and Section~3.

In a similar fashion one can define ``triple``~ $\beta$-Grothendieck 
polynomials of the classical type $W=A,B,C,D$:

$$({\mathfrak{G}}^{A_{n-1}}(-Y_n))^{-1}~{B^{W}(Z_m)}~{\mathfrak{G}}^{A_{n-1}}(-X_n)=
\sum_{w \in W}~{\mathfrak{G}}^{W}_{w}(X_n,Y_n,Z_m)~u_{w},$$ 
where now the left and the right parts are treated in the Id-Coxeter algebra 
$Id_{\beta}(W)$.

 An algebra-combinatorial approach is used as the basic tool in the present 
paper gives rise naturally to the study of the generating functions for the 
double and triple Schubert and $\beta$-Grothendieck polynomials are 
introduced  in the present notes, \underline{but} in more wider class of 
algebras such as 
the Hecke and Temperley--Lieb algebras of classical types, and the plactic 
and reduced plactic algebras of classical type. Recall that the plactic 
algebra of classical type $W$ is the quotient of the unital free associative 
algebra over  $\Q$ of rang $n:=r(W)$ by the two-sided ideal generated by the 
$W$-plactic (or $W$-Knuth--Kra\'{s}kiewicz) relations has been described and
studied in \cite{Le}.    

\begin{prob} ${}$

$({\bf A})$~~~To extend  results concerning the plactic algebra and plactic 
polynomials of type $A_{n-1}$ obtained in  \cite{Ki} to the case of the 
plactic algebras corresponding to plactic monoid of type $ W:= B_n$, $C_n$
 and  $D_n$, introduced in \cite{Le}.~~In particular, describe the MacNeille 
completion ${\cal{MN}}(W)$  of the the Bruhat graph (poset) associated with 
the Weyl groups of classical type, as well as describe the decomposition of 
the $W$-Cauchy kernel in the reduced $W$-plactic algebra.

$({\bf B})$~~~Find a geometric interpretation of plactic Schubert and 
Grothendieck polynomials of classical types. Does the MacNeille complition 
${\cal{MN}}(W)$ can be realized as a convolution algebra of a certain 
nonsingular algebraic variety ?

\end{prob}
A few words about the history of \underline{problems} considered in the 
present paper in order. The algebraic and combinatorial theory of single and 
double Schubert polynomials of type $A$ was initiated and studied 
comprehensively by 
A. Lascoux and M.-P. Sch\"{u}tzenberger in the middle of 80's of the last 
century.  We refer the reader to the nice written books \cite{Mac} and 
\cite{Man} for detailed exposition of this subject. The general description 
of the cohomology
 and  equivariant cohomology rings, K-theory and equivariant K-theory of 
generalized flag varieties corresponding to a symmetrizable Kac--Moody group 
was created by B. Kostant and S. Kumar. Details can be found in the book 
 \cite{Ku} \\

A bit of history concerning the present notes. ~This paper 
(as well as \cite{Ki}) is an update version of 
my notes written for Course `` Schubert Calculus''~has been delivered at the 
Graduate School of Mathematical Sciences, University of Tokyo (1995/96), and 
 at the Graduate School of Mathematics, Nagoya University (1998/99).\\   

Final remark, in \cite{IMN} the polynomials $\s_{w}^{W}(X,Y,Z)$ have been 
rediscovered using a geometrical approach, see also \cite{AF}.

\section{Basic definitions}
\subsection{Weyl groups of classical types}
\subsubsection{The symmetric group}
The symmetric group ${\mathbb S}_n,$~$n \ge 1,$ is the group of all 
permutations of the set $[1,n]:= \{ 1,2,\ldots, n \}.$  As is customary, we 
will identify a permutation $w \in  {\mathbb S}_n$ with its image, i.e. with 
the sequence $(w(1), w(2), \ldots, w(n)).$ Sometimes we will write $w_i$ 
instead of $w(i),$ and $w_1 \ldots w_n$ instead of sequence 
$(w(1), w(2), \ldots, w(n)).$

For $i=1,\ldots,n-1$  let $s_i$ denote the transposition that interchanges $i$ 
and $i+1,$ and fixes all other numbers in $[1,n].$ It is well-known that the
elements $s_1,\ldots, s_{n-1}$ generate the symmetric group ${\mathbb S}_n$ and
 satisfy the following relations

$(1)$~~~$s_i^2=1;$

$(2)$~~$s_i~s_j=s_j~s_i,$~~if~~$|i-j| \ge 2;$

$(3)$~~$s_i~s_{i+1}~s_{i}=s_{i+1}~s_{i}~s_{i+1}$~for~$i=1,\ldots,n-2.$

For a permutation $w \in {\mathbb S}_n$ let's  denote by $D(w)$ the {\it 
diagram} of the permutation $w,$ see e.g. \cite{Mac}, i.e.
$$ (i,j) \in D(w) \Longleftrightarrow  i< w^{-1}(j) ~~and~~j < w(i).$$
It is well-know that $l(w)= |D(w)|,$ where $l(w)$ denotes the length of a 
permutation $w,$ i.e. the minimal number of generators whose product is $w.$  

\subsubsection{The hyperochtahedral group}
The hyperochtahedral group $B_n:=W(B_n),$~$ ~n \ge 2,$ is the group of 
symmetries of the $n$-dimensional cube. As an abstract group it can be given 
by the set of generators $s_{0},s_{1},\ldots,s_{n-1}$ satisfying relations

$(1)$~~$s_i^2=1,$ ~~if ~~$i=0,1,2,\ldots, n-1;$

$(2)$~~$s_i~s_j=s_j~s_i$, ~~if ~$|i-j| \ge 2;$

$(3)$~~$s_i~s_{i+1}~s_i=s_{i+1}~s_{i}~s_{i+1},$ ~~if~$i=1,\ldots,n-2;$

$(4)$~~$s_{0}~s_1~s_{0}~s_1=s_1~s_{0}~s_1~s_{0}.$

The elements of $B_n$ can be thought of as {\it signed permutations}: a 
generator $s_i,$ $i > 0,$ interchanges entries in the $i'$th and 
$(i+1)'$st positions, and the generator $s_{0}$ changes the sign of the first 
entry. As in any Coxeter group, the {\it length l(w)} of an element $w$ is the minimal number of generators whose product is $w.$  Such a factorization of 
minimal length, or the corresponding sequence of indices, is called a {\it 
reduced decomposition} of $w.$ As is customary, we will write ${\bar i}$ 
instead of $-i.$ For example, the action of the generator $s_0$  looks like
~$s_{0}(1 2\ldots n)={\bar 1} 2 \ldots n,$ and ${\bar {\bar a}}= a.$ For any 
sign permutation $w=w_1w_2 \ldots w_n $ let $ {\bar w}$ denotes the sign 
permutation ${\bar w_1} {\bar w_2} \ldots {\bar w_n}.$ It is clear that if 
$w \in {\mathbb S}_n \subset B_n,$ then $l(w)+l({\bar w})=n^2.$  

\subsubsection{The group $W(D_n)$}
The group $D_n:=W(D_n)$ is a subgroup of elements $w \in W(B_n)$ which make an 
even number of sign changes. The standard generators for this group are $s_i,$ 
$i=1, \ldots, n-1$ and $s_{{\hat 1}},$  subject to the set of relations
$$\left\{\matrix{ 
s_i^2=1, & if & i={\hat 1},1,2,\ldots, n-1; \cr
s_i~s_j=s_j~s_i, & if & |i-j| \ge 2; \cr
s_i~s_{{\hat 1}}=s_{{\hat 1}}~s_i, & if & i \not= 2; \cr
s_i~s_{i+1}~s_i=s_{i+1}~s_{i}~s_{i+1}, & if & i=1,\ldots,n-2; \cr
s_2~s_{{\hat 1}}~s_2=s_{{\hat 1}}~s_2~s_{{\hat 1}}. \cr} \right. $$
The elements of $D_n$ can be thought of as {\it even signed permutations}: a 
generator $s_i,$ $i > 0,$ interchanges variables $x_i$ and $x_{i+1}$, and the 
generator $s_{{\hat 1}}$ replaces $x_1$ with $-x_2$ and $x_2$ with $-x_1.$

\subsection{Divided difference and isobaric divided difference operators}

The hyperoctahedral group $B_n$ acts on the ring of polynomials 
$ P_n:=\Q[x_1, \ldots, x_n]$ in the natural way. Namely, $s_i$ interchanges 
$x_i$ and $x_{i+1},$  for $i=1,\ldots,n-1$, and the special generator $s_{0}$ 
acts by
$$s_{0}f(x_1,x_2,\ldots,x_n ) = f(-x_1,x_2,\ldots,x_n ) .$$  
The divided difference operator $\partial_{i}$ for $i=1,\ldots,n-1,$ acts on 
the ring of polynomials $P_n$ by
$$\partial_i f(x_1,\ldots,x_n)={f -s_i(f) \over x_i -x_{i+1}} ,$$
and the $B$ type divided difference operator $\partial_{0}:=\partial_{0}^{B}$ 
acts on the ring of polynomials $P_n$ by
$$ \partial_{0} f(x_1,x_2,\ldots,x_n)= 
{f(x_1,x_2,\ldots,x_n)- f(-x_1,x_2,\ldots,x_n) \over x_1}.$$
We consider also $C$ and $D$ types divided difference operators 
$\partial_{0}^{C}$ and $\partial_{{\hat 1}}^{D}$ which act on the ring of 
polynomials by ~~~$ \partial_{0}^{C} (f):={1/2}~\partial_{0}^{B} (f)$~~and
$$\partial_{{\hat 1}}^{D} f(x_1,x_2,\ldots,x_n)= 
{f(x_1,x_2,\ldots,x_n)- f(-x_2,-x_1,\ldots,x_n)\over x_1+x_2}.$$

Finally, we define {\it isobaric divided difference operators} 
$\pi_{\alpha}^{G}$  for each simple root $\alpha$ of the corresponding Lie 
group $G$ of classical type. Namely, let $\beta$ be a parameter, define
 ($ 1 \le i \le n-1$)
$$ \pi_{i}^{A}(f)= \partial_{i}((1+\beta~x_{i+1}) f)=\pi_{i}^{C}(f),~~
  \pi_{i}^{B}(f)= \partial_{i}((1+{\beta \over 2}~x_{i+1})~f)=
\pi_{i}^{D}(f),$$
$$\pi_{0}^{B}(f)=\partial_{0}^{B}((1-{\beta \over 2}~x_1) f),~~
\pi_{0}^{C}(f)=\partial_{0}^{C}((1-{\beta}~x_1) f),$$
$$\pi_{{\hat 1}}^{D}(f)= 
\partial_{{\hat 1}}^{D}((1-{\beta \over 2}~x_1)(1-{\beta \over 2}~x_2) f).$$
Note that $\pi_{i}^2=- \beta \pi_{i}$ for types $A$ and $C;$~~
$\pi_{i}^2=- {\beta \over 2} \pi_{i}$ for types $B$ and $D;$~~
$\pi_{0}^2=- {\beta} \pi_{0}$ ~for type $B;$~~
$\pi_{0}^2=- 2 {\beta} \pi_{0}$ ~for type $C,$~and~
$\pi_{{\hat 1}}^2=- {\beta} \pi_{{\hat 1}}$ for type $D.$

\subsection{Nil--Coxeter and Id--Coxeter algebras of classical types}
\subsubsection{Nil--Coxeter  and Id--Coxeter algebras of type $A$} 
$\bullet$~~~Let $NC_n$ denotes the Nil--Coxeter algebra type $A_{n-1}.$ 
Recall that $NC_n$ is an associative algebra generated over $\Z$ by the set of 
generators $\{u_1, \ldots,u_{n-1} \}$ subject to the set of relations

$(a)$~~~$u_i^{2}=0$ ~for $i=1, \ldots, n-1,$

$(b)$~~$u_i~u_j=u_j~u_i,$ ~if ~$1 \le i,j \le n-1$~and~$|i-j| \ge 2,$ 

$(c)$~~$u_i~u_{i+1}~u_i=u_{i+1}~u_{i}~u_{i+1},$~if~$i=1, \ldots,n-2.$

It is well-known that $dim NC_n = n !$ and the elements 
$\{ u_{w}, w \in {\mathbb S}_n \}$ form a $\Z$-linear basis in the algebra  
$NC_n,$ where by definition we set $u_{w}=u_{a_1} \ldots u_{a_l}$ for any 
reduced decomposition $w=s_{a_1} \ldots u_{a_l}$  of $w \in {\mathbb S}_n$ 
chosen.

$\bullet$~~~Let $\beta$ be a parameter. The Id--Coxeter algebra of type $A,$ 
denoted by $Id(A_{n-1}):=Id_{\beta}(A_{n-1}),$ is 
an associative algebra generated over $\Z[\beta]$ by the set of generators 
$\{u_1, \ldots,u_{n-1} \}$ subject to the set of relations $(b)$ and $(c)$  
from the definition of the algebra $NC_n,$ and the relations $u_i^2 =\beta~u_i$
~for $i=1,\ldots,n-1,$ instead of that $(a).$ It is well-known that the 
elements $\{ u_{w}, w \in {\mathbb S}_n \}$ form a $\Z[\beta]$-linear basis 
in the algebra  $Id_{\beta}(A_{n-1}).$ 

\subsubsection{Nil--Coxeter and Id--Coxeter algebras of type $B$} 

$\bullet$~~~Let $Nil(B_n)$ denotes the nil--Coxeter algebra of type $B.$ 
Recall that 
$Nil(B_n)$ is an associative algebra generated over $\Z$ by the set of 
generators $\{u_0,u_1, \ldots,u_{n-1} \}$ subject the set of relations

$(a)$~~~$u_i^{2}=0$ ~for $i=0,1, \ldots, n-1,$

$(b)$~~$u_i~u_j=u_j~u_i,$ ~if ~$1 \le i,j \le n-1$~and~$|i-j| \ge 2,$ 

$(c)$~~$u_i~u_{i+1}~u_i=u_{i+1}~u_{i}~u_{i+1},$~if~$i=1, \ldots,n-2,$

$(d)$~~$u_{0}~u_1~u_{0}~u_{1}=~u_1~u_{0}~u_{1}~u_{0},$~~$u_{0}~u_i=u_i~u_{0}$ 
for $i=2,\ldots,n-1.$

It is well-known that $dim Nil(B_n) = 2^{n}~n !$ and the elements 
$\{ u_{w}, w \in B_n \}$ form a $\Z$-linear basis in the algebra  $Nil(B_n),$ 
where by definition we set $u_{w}=u_{a_1} \ldots u_{a_l}$ for any reduced 
decomposition $w=s_{a_1} \ldots u_{a_l}$  of $w \in W(B_n)$ chosen.

$\bullet$~~~Let $\beta$ be a parameter. The Id--Coxeter algebra of type $B,$ 
denoted by $Id(B_n):=I_{\beta}(B_n),$ is an associative algebra generated 
over the ring $\Q[\beta]$ by the set of generators 
$\{u_{0},u_1, \ldots,u_{n-1} \}$ 
subject to the set of relations $(b),$ $(c)$ and $(d)$ from the definition of
 the algebra $Nil(B_n),$ and the relations $u_i^2 =\beta~u_i$~for $i=0,1,
\ldots,n-1,$ instead of that $(a).$ It is 
well-known that the elements $\{ u_{w}, w \in W(B_n) \}$ form a 
$\Z[\beta]$-linear basis in the algebra  $Id_{\beta}(B_n).$ 

Let $x_1,\ldots,x_n$ be a set of variables which assumed to be commute with 
all generators $u_0, \ldots, u_{n-1}.$ Define deformed addition $x +_{\beta}~y=
x+y+\beta~x~y,$ so that $x-_{\beta}~y= (x-y)/(1+\beta~y).$

Follow \cite{FK3}, define 
$$h_i(x):= 1+x~u_i,~~~~for~~~i=1, \ldots, n-1,~~~
h_{0}(x):= h_{0}^{B}(x)=1+x~u_0.$$
Define also 
$$h_{0}^{C}(x)=1+2 x~u_{0}~~~and ~~~
h_{{\hat 1}}^{D}(x)=1+x~u_{{\hat 1}}:=h_{{\hat 1}}(x).$$
\begin{lem} (Cf  \cite{FS},\cite{FK3})~~The elements $h_{i}(x)$ satisfy the 
following relations

$(1)$~~~$h_{i}(x)~h_{j}(y)=h_{j}(y)~h_{i}(x),$~if
~$1 \le i,j \le n-1$~and~$|i-j| \ge 2,$

$(2)$~~$h_{0}(x)~h_{i}(y)= h_{i}(y)~h_{0}(x),$~if $i=2, \ldots, n-1,$

$(2)$~~$h_i(x)~h_i(y)=h_{i}(x+y)$ ~~~in the algebra $Nil(B_n),$

$(2a)$~~$h_i(x)~h_i(y)=h_{i}(x+y+\beta ~x~y)=h_i(x+_{\beta}y)$~~~
in the algebra $Id_{\beta}(B_n),$

$(3)$~~( Yang--Baxter equation of type $A$ in the algebra $Nil(B_n)$)
$$h_{i}(x)~h_{i+1}(x+y)~h_{i}(y)= h_{i+1}(y)~h_{i}(x+y)~h_{i+1}(x),
~~i=1,\ldots,n-2,$$

$(3a)$~~( Yang--Baxter equation of type $A$ in the algebra 
$Id_{\beta}(B_n)$ )
$$h_{i}(x)~h_{i+1}(x+_{\beta}y)~h_{i}(y)= 
h_{i+1}(y)~h_{i}(x+_{\beta}y)~h_{i+1}(x),
~~i=1,\ldots,n-2,$$

$(4)$~~( Yang--Baxter equation of type $B$ in the algebra $Nil(B_n)$)
$$h_{0}(y)~h_{1}(x+y)~h_{0}(x)~h_{1}(x-y)=
h_{1}(x-y)~h_{0}(x)~h_{1}(x+y)~h_{0}(y).$$

$(4a)$~~( Yang--Baxter equation of type $B$ in the algebra 
$Id_{\beta}(B_n)$)
$$h_{0}(y)~h_{1}(x+_{\beta}y)~h_{0}(x)~h_{1}(x-_{\beta}y)=
h_{1}(x-_{\beta}y)~h_{0}(x)~h_{1}(x+_{\beta}y)~h_{0}(y).$$
\end{lem}
Let us introduce in the algebra ~$Id_{\beta}(B_n)$~the elements~(cf \cite{FS},
 \cite{FK2}):
$$ {\cal A}(x):=A_{i}(\phi_{\beta}(x)):= A_{i}^{(n)}(\phi_{\beta}(x))=
\prod_{a =n-1}^{i}~h_{a}(\phi_{\beta}(x)),$$
$${\cal B}(x):=B(\phi_{\beta}(x)):= \prod_{a =n-1}^{1}~h_{a}(\phi_{\beta}(x))
~h_{0}(\phi_{\beta}(x))~\prod_{a =1}^{n-1}~h_{a}(\phi_{\beta}(x)),$$
$${\cal C}(x):=C(\phi_{2\beta}(x)):= \prod_{a =n-1}^{1}~h_{a}(\phi_{2\beta}(x))
~h_{0}(\phi_{2\beta}(2x))~\prod_{a =1}^{n-1}~h_{a}(\phi_{2\beta}(x)),$$
where $\phi_{\beta}(x):=x / (1-{\beta \over 2}x).$

In the nil--Coxeter algebra $Nil(B_n)$ 
these elements can be written in the form
$$B_n(x):= A_{1}^{(n)}(x)~h_{0}(x)~A_{1}^{(n)}(-x)^{-1},$$
$$ C(x):=C_n(x)= A_{1}^{(n)}(x)~h_{0}(2x)~A_{1}^{(n)}(-x)^{-1}.$$
\begin{lem} ${}$ (\cite{FS}, \cite{FK3}) ~~One has

$(1)$~~$A_{i}(x)~A_{i}(y)=A_{i}(y)~A_{i}(x),$~~$A_{i}(x)=A_{i+1}(x)~h_{i}(x).$

$(2)$~~$B(x)~B(y)=B(y)~B(x),$~~~$C(x)~C(y)=C(y)~C(x)$

in the both algebras $Nil(B_n)$ and $Id_{\beta}(B_n).$~~Therefore,

$(2a)$~~${\cal B}(x)~{\cal B}(y)={\cal B}(y)~{\cal B}(x),$~~~
${\cal C}(x)~{\cal C}(y)={\cal C}(y)~{\cal C}(x)$

in the algebra $Id_{\beta}(B_n).$

$(3)$~~$B(x)~B(-x)=1,$~~~$C(x)~C(-x)=1$ in the Nil--Coxeter algebra 
$Nil(B_n).$ 

$(3a)$~~${\cal B}(x)~{\cal B}(-x)=1,$~~~${\cal C}(x)~{\cal C}(-x)=1$ in the 
algebra $Id_{\beta}(B_n).$
\end{lem} 
Finally,  let us consider the following expressions in the algebra $Nil(B_n)$ 
\begin{equation}
H(Z_m)=B(z_1)~B(z_2) \cdots B(z_m)= \sum_{w \in W(B_n)}~F_{w}(Z_m)~u_{w},
\end{equation} 
\begin{equation}
{\HH}(t_1, t_2,\ldots, t_m) = {\sqrt  {B(t_1)~B(t_2) \cdots B(t_m)} },
\end{equation}
and those in the algebra $Id_{\beta}(B_n)$
\begin{equation}
{\cal H}(Z_m)={\cal B}(z_1)~{\cal B}(z_2) \cdots {\cal B}(z_m)= 
\sum_{w \in W(B_n)}~{\cal F}_{w}(Z_m)~u_{w},
\end{equation} 
\begin{equation}
{\H}(t_1, t_2,\ldots, t_m) = 
{\sqrt  { {\cal B}(t_1)~{\cal B}(t_2) \cdots {\cal B}(t_m) }} .
\end{equation}
It follows from Lemma~2.2 that the  $H(Z_m)$ and ${\cal H}(Z_m)$ as well as 
${\HH}(t_1,\cdots,t_m)$ and  ${\H}(t_1, t_2,\ldots, t_m),$  
are symmetric functions of the variables $z_1, \ldots,z_m$ and 
$t_1, \ldots, t_m$ respectively.
\begin{lem} ${}$

$(1)$ For any $w \in W(B_n)$ the polynomials $F_{w}(Z_m)$ and 
${\cal F}_{w}(Z_m)$ are  supersymmetric 
functions of the variables $Z_m=(z_1, \ldots,z_m),$ i.e. $F_{w}(Z_m)$ and 
${\cal F}_{w}(Z_m)$ are polynomials  of the odd power sums 
$p_1(Z_m),p_3(Z_m), \ldots.$

$(2)$~(\cite{BH}, \cite{Lam})  For any $w \in W(B_n),$ polynomial 
~$F_{w}(Z_m)$ is a linear combination with non-negative integer coefficients 
of Schur $P$-functions.

$(3)$~(\cite{FK3})~~Assume that the variables $z_1,z_2, \ldots $  and 
$t_1,t_2, \ldots $ are related by
\begin{equation}
{p_k(t_1,t_2, \ldots) \over 2}= p_k(z_1,z_2, \ldots), ~~k=1,3,5, \ldots. 
\end{equation} 
Then ~~${\HH}(t_1,t_2,\ldots )=H(z_1,z_2, \ldots).$
\end{lem}
\begin{ex} We display the polynomials $F_{w}(Z_m)$ for $n=2,m=4. $\\
$F_{id}(Z_4)=1,$ \\
$F_{u_{0}}(Z_4)=z_1+z_2+z_3+z_4,$ \\
$F_{u_{1}}(Z_4)=2(z_1+z_2+z_3+z_4),$ \\
$F_{u_{01}}(Z_4)=F_{u_{10}}(Z_4)=(z_1+z_2+z_3+z_4)^2,$ \\
$F_{u_{010}}(Z_4)=z_1~z_2(z_1+z_2)+(z_1+z_2)(z_3+z_4)(z_1+z_2+z_3+z_4)+
z_3~z_4(z_3+z_4),$ \\
$F_{u_{101}}(Z_4)=(z_1+z_2)(z_1^2+z_1~z_2+z_2^2)+z_3~z_4(z_3^2+z_3~z_4+
z_4^2)+$ \\ 
$2(z_1+z_2)(z_3+z_4)(z_1+z_2+z_3+z_4),$ \\
$F_{u_{0101}}(Z_4)= (z_1+z_2+z_3+z_4)~F_{u_{010}}(Z_4).$
\end{ex}
Note that $F_{u_{010}}(Z_4)= s_{(2,1)}(Z_4).$

\subsubsection{Nil--Coxeter and Id--Coxeter algebras of type $D$} 
$\bullet$~~Denote by $Nil(D_n)$ the nil--Coxeter algebra type $D.$ Recall 
that $Nil(D_n)$ 
is an associative algebra generated over $\Q$ by the set of elements 
$\{ u_{{\hat 1}}, u_{1},u_{2},\ldots,u_{n-1} \}$ subject to the set of 
relations 

$(a)$~~$u_{{\hat 1}}^2=0,$~~~$u_{i}^2=0$ ~for $i=1,2,\ldots,n-1;$

$(b)$~~$u_i~u_j=u_j~u_i,$ ~~if ~$|i-j| \ge 2;$

$(c)$~~$u_{i}~u_{{\hat 1}}=u_{{\hat 1}}~u_{i},$~~if~$i \not= 2;$~~
$u_{{\hat 1}}~u_{2}~u_{{\hat 1}}=u_{2}~u_{{\hat 1}}~u_{2},$

$(d)$~~$u_{i}~u_{i+1}~u_{i}=u_{i+1}~u_{i}~u_{i+1}$~for $i=1,2,\ldots,n-2.$

It is well-known that the elements $\{ u_{w}, w \in W(D_n) \}$ form 
a $\Z[\beta]$-linear basis in the algebra  $Nil(D_n).$ 

$\bullet$~~Let $\beta$ be a parameter. The Id--Coxeter algebra $Id(D_n):=
Id_{\beta}(D_n)$ is 
an associative algebra generated over $\Q[\beta]$ by the set of generators 
$\{u_{{\hat 1}},u_1, \ldots,u_{n-1} \}$ subject to the set of relations 
$(b),$ $(c)$ and 
$(d)$ from the definition of the algebra $Nil(D_n),$ and the relations 
$u_i^2 =\beta~u_i$~for $i={\hat 1},1,\ldots,n-1,$ instead of that $(a).$ 
It is well-known that the elements $\{ u_{w}, w \in W(D_n) \}$ form 
a $\Z[\beta]$-linear basis in the algebra  $Id_{\beta}(D_n).$ 

Define ~~$D(x):=D_n(x)=h_{n-1}(x) \cdots h_{1}(x)~h_{{\hat 1}}(x)~h_{2}(x) 
\cdots h_{n-1}(x).$

Recall that $h_{{\hat 1}}(x)=1+x~u_{{\hat 1}}.$

First we study properties of the elements $D(x)$ in the Nil--Coxeter algebra 
$Nil(D_n).$ It is easy to see that
$$D(x)=A_{1}(x)~h_{{\hat 1}}(x)~A_{2}(-x)^{-1}=
A_{2}(x)~h_{{\hat 1}}(x)~A_{1}(-x)^{-1}.$$
\begin{lem}~~~$D(x)~D(y)=D(y)~D(x),$~~$D(x)~D(-x)=1.$
\end{lem}
{\bf Proof.} ~~~One has ~~~~$D(x)~D(y)=$

$ A_1(x)~h_{{\hat 1}}(x)~A_{2}(-x)^{-1}~A_{2}(y)~h_{{\hat 1}}(y)
~A_{1}(y)^{-1}=$ \\
$A_{1}(x)~A_{1}(y)~A_{1}(y)^{-1}~h_{{\hat 1}}(x)~A_{2}(y)~A_{2}(-x)^{-1}~
h_{{\hat 1}}(y)~A_{1}(-x)~A_{1}(-x)^{-1}~A_{1}(-y)^{-1}=$ \\
$A_{1}(x)~A_{1}(y)~h_{1}(-y)~h_{2}(-y)~h_{{\hat 1}}(x)~A_{3}(y)^{-1}A_{2}(y)~
A_{2}(-x)^{-1}~A_{3}(-x)$ \\
$h_{{\hat 1}}(y)~h_{2}(-x)~~h_{1}(-x)~A_{1}(-x)^{-1}~A_{1}(-y)^{-1}=$ \\
$A_{1}(x)~A_{1}(y)~h_{1}(-y)~h_{2}(-y)~h_{{\hat 1}}(x)~h_{2}(x+y)~
h_{{\hat 1}}(y)~h_{2}(-x)~h_{1}(-x)$ \\
$~A_{1}(-x)^{-1}~A_{1}(-y)^{-1}= A_{1}(x)~A_{1}(y)~h_{{\hat 1}}(x+y)
~h_{1}(-x-y)~A_{1}(-x)^{-1}~A_{1}(-y)^{-1}.$

The final expression is symmetric with respect to $x$ and $y,$ therefore the 
elements $D(x)$ and $D(y)$ commute to one another.

Note that to deduce the final equality we have used the Yang--Baxter relation
$$h_{2}(x)~h_{{\hat 1}}(x+y)~h_{2}(y)=
h_{{\hat 1}}(x)~h_{2}(x+y)~h_{{\hat 1}}(y).$$
$\qed$
\begin{lem} The elements $D(x)$ and $D(y)$ commute in the algebra $Id(D_n)$ and
~~~$D(\phi_{\beta}(x))~D(\phi_{\beta}(-x))=1,$ 
where $\phi_{\beta}(x)=x/(1-{\beta \over 2}~(x).$ 
\end{lem}
{\bf Proof.} It is clear that the element $D_2(x)=h_1(x)~h_{{\hat 1}}(x)$ 
commutes with $D_2(y).$ The next case is $n=3.$ We have \\
$D_3(x)~D_3(y)=h_2(x)~h_1(x)~h_{{\hat 1}}(x)~h_2(x+_\beta y)~h_{{\hat 1}}(y)~h_1(y)~h_2(y)=$\\
$h_2(y)~h_2(x-_{\beta}y)~h_1(x)~h_2(y)~h_{{\hat 1}}(x+_{\beta} y)
~h_2(x)~h_{1}(y)~h_{2}(y-_{\beta} x)~h_2(x)=$\\
$h_2(y)~h_1(y)~h_{{\hat 1}}(x+_{\beta} y)~h_1(x)~h_2(x)=D_3(x)~D_3(y).$\\
 Now we make use an induction. We have ~~~~$D_{n+1}(x)~D_{n+1}(y)=$ \\
$h_{n}(x)~h_{n-1}(x)~D_{n-1}(x)~h_{n-1}(x)~h_{n}(x+_{\beta} y)~h_{n-1}(y)
~D_{n-1}(y)~h_{n-1}(y)~h_{n}(y)=$\\
$h_{n}(x)~h_{n-1}(x)~h_{n}(y)~
D_{n-1}(x)~h_{n-1}(x+_{\beta} y)~D_{n-1}(y)~h_{n}(x)~h_{n-1}(y)~h_{n}(y)=$\\ 
$h_{n}(x)~h_{n-1}(x)~h_{n}(y)~h_{n-1}([-x]_{\beta})
D_{n}(x)~D_{n}(y)~h_{n-1}([-y]_{\beta})~h_{n}(x)~h_{n-1}(y)~h_{n}(y)=$\\ 
$h_{n}(x)~h_{n-1}(x)~h_{n}(y)~h_{n-1}([-x]_{\beta})
D_{n}(y)~D_{n}(x)~h_{n-1}([-y]_{\beta})~h_{n}(x)~h_{n-1}(y)~h_{n}(y)=$\\ 
$h_{n}(x)~\biggl(h_{n-1}(x)~h_{n}(y)~h_{n-1}(y-_{\beta} x) \biggr)
D_{n-1}(y)~h_{n-1}(x+_{\beta} y)~D_{n-1}(x)$ \\
$\biggl( h_{n-1}(x-_{\beta} y)~h_{n}(x)~h_{n-1}(y) \biggr)~h_{n}(y)=$\\ 
$h_{n}(y)~h_{n-1}(y)~h_{n}(x)~D_{n-1}(y)~h_{n-1}(x+_{\beta} y)~D_{n-1}(x)~
h_{n}(y)~h_{n-1}(x)~h_{n}(x)=$\\ 
$D_{n+1}(y)~D_{n+1}(x).$  \\
The second statement follows from the identity 
$$\phi_{\beta}(x)+_{\beta} \phi_{\beta}(-x)=0.$$
$\qed$

Let us set ${\cal D}(x):={\cal D}_n(x)=D(\phi_{\beta}(x)).$ It follows from 
Lemma~2.5 that in the algebra $Id_{\beta}(D_n)$ the elements ${\cal D}(x)$ 
and ${\cal D}(y)$  commute, and ${\cal D}(x)~{\cal D}(-x)=1.$

\section{Schubert and Grothendieck polynomials}
\subsection{Schubert and Grothendieck expressions}
Let us recall, cf \cite{FS}, \cite{FK2},  the definition of certain elements 
in the Nil--Coxeter and Id--Coxeter algebras which will be used to define 
Schubert and Grothendieck expressions. 
$$ A_{i}(x):= A_{i}^{(n)}(x)=\prod_{a =n-1}^{i}~h_{a}(x),~~
B(x):=B_n(x)= \prod_{a =n-1}^{1}~h_{a}(x)~h_{0}(x)~
\prod_{a =1}^{n-1}~h_{a}(x) ,$$
$$ C(x):=C_n(x)= A_{1}^{(n)}(x)~h_{0}(2x)~A_{1}^{(n)}(-x)^{-1}.$$
 
\begin{dfn}~~~(The case of Nil--Coxeter algebras) ${}$

$(A)$ ~(\cite{FS})~~The type $A$ Schubert expression is
$$ \s^{A_{n-1}}(X_{n}):=A_{1}^{(n)}(x_1)~A_{2}^{(n)}(x_2)~
\cdots A_{n-1}^{(n)}(x_{n-1}) .$$
$(B)$~~The type $B$ Schubert expression  of the first kind is
$$ \s^{B_{n}}(X_n):= B_n(x_1)~B_n(x_2) \cdots B_n(x_n)~\s^{A_{n-1}}(- X_n).$$
$(Ba)$~~The type $B$ Schubert expression of the second kind is 
$$ \b (X_n):= \sqrt {B_n(x_1)~B_n(x_2) \cdots B_n(x_n)}~
\s^{A_{n-1}}(- X_n),$$
$(C)$~~The type $C$ Schubert expression  of the first kind is
$$ \s^{C_{n}}(X_n):= C_n(x_1)~C_n(x_2) \cdots C_n(x_n)~\s^{A_{n-1}}(- X_n).$$
$(Ca)$~~The type $C$ Schubert expression of the second kind is 
$$ \c (X_n):= \sqrt {C_n(x_1)~C_n(x_2) \cdots C_n(x_n)}~
\s^{A_{n-1}}(- X_n),$$
$(D)$~~The type $D$ Schubert expression of the first kind is 
$$ \s^{D_{n}}(X_n):= D_n(x_1) \cdots D_n(x_{n-1})~A_2(-x_1)~A_3(-x_2)
 \cdots A_{n-1}(-x_{n-2}),$$
$(Da)$~~The type $D$ Schubert expression of the second kind is 
$$ \D (X_n):= \sqrt {D_n(x_1) \cdots D_n(x_{n})}~\s^{A_{n-1}}(-X_n).$$
\end{dfn}
In a similar fashion one can define $\s^{A}(X_{\infty})$,~
$\s^{B}(X_{\infty}),$~$\b(X_{\infty}),$~$\s^{C}(X_{\infty}),$ \\
$\c(X_{\infty}),$~$\s^{D}(X_{\infty})$ and $\D (X_{\infty}).$ 
\begin{dfn} (The case of $Id$-Coxeter algebras)

$(A)$ ~~The type $A$ Grothendieck  expression is
$$ \g^{A_{n-1}}(X_{n}):=A_{1}^{(n)}(x_1)~A_{2}^{(n)}(x_2)~
\cdots A_{n-1}^{(n)}(x_{n-1}) .$$
$(B)$~~The type $B$ Grothendieck  expression  of the first kind is
$$ \g^{B_{n}}(X_n):= {\cal B}_n(x_1)~{\cal B}_n(x_2) \cdots {\cal B}_n(x_n)
~\s^{A_{n-1}}(\phi_{\beta}(- X_n)).$$
$(Ba)$~~The type $B$ Grothendieck expression of the second kind is 
$$ \g_{B_n} (X_n):= \sqrt {{\cal B}_n(x_1)~{\cal B}_n(x_2) \cdots 
{\cal B}_n(x_n)}~\s^{A_{n-1}}(\phi_{\beta}(- X_n)),$$
$(C)$~~The type $C$ Grothendieck expression  of the first kind is
$$ \g^{C_{n}}(X_n):= {\cal C}_n(x_1)~{\cal C}_n(x_2) \cdots {\cal C}_n(x_n)~
\s^{A_{n-1}}(\phi_{2 \beta}(- X_n)).$$
$(Ca)$~~The type $C$ Grothendieck  expression of the second kind is 
$$ \g_{C} (X_n):= \sqrt {{\cal C}_n(x_1)~{\cal C}_n(x_2) \cdots 
{\cal C}_n(x_n)}~\s^{A_{n-1}}(\phi_{2\beta}(- X_n)),$$
$(D)$~~The type $D$ Grothendieck expression of the first kind is~~~ 
$ \g^{D_{n}}(X_n):=$ 
$${\cal D}_n(x_1) \cdots {\cal D}_n(x_{n-1})~
A_2(\phi_{\beta}(-x_1))~A_3(\phi_{\beta}(-x_2))\cdots 
A_{n-1}(\phi_{\beta}(-x_{n-2}),$$
$(Da)$~~The type $D$ Grothendieck expression of the second kind is 
$$ \g_{D_{n}} (X_n):= \sqrt {{\cal D}_n(x_1) \cdots {\cal D}_n(x_{n}) }
~\s^{A_{n-1}}(\phi_{\beta}(-X_n)).$$
\end{dfn}
\begin{rem}{\rm Since the Nil--Coxeter algebras in question are 
finite dimensional ( in fact, the both Nil--Coxeter and Id--Coxeter algebras 
do not have elements of degree $ > n$ ), the Schubert and Grothendieck 
expressions of the second kind are {\it polynomials}.  
}
\end{rem}

\begin{dfn} (\cite{BH}, \cite{FK3}) ${}$ 

$(A)$ ~~Schubert polynomials of type $A$ are defined via the decomposition
$$ \s^{A_{n-1}}(X_n)= \sum_{w \in {\mathbb S}_n}~ \s_{w}(X_n)~u_{w},$$
$(B)$~~The type $B$ Schubert polynomials of the first kind are defined via the 
decomposition
$$ \s^{B_{n}}(X_n)= \sum_{w \in W(B_n)}~\s_{w}^{B_{n}}(X_n)~u_{w},$$
$(Ba)$~~The type $B$ Schubert polynomials of the second kind are defined via 
the decomposition
$$ \b(X_n)= \sum_{w \in W(B_n)}~\b_{w}(X_n)~u_{w},$$
$(C)$~~The type $C$ Schubert polynomials of the first kind are defined via the 
decomposition
$$ \s^{C_{n}}(X_n)= \sum_{w \in W(C_n)}~\s_{w}^{C_{n}}(X_n)~u_{w},$$
$(Ca)$~~The type $C$ Schubert polynomials of the second kind are defined via 
the decomposition
$$ \c(X_n)= \sum_{w \in W(C_n)}~\c_{w}(X_n)~u_{w},$$
$(D)$~~The type $D$ Schubert polynomials of the first kind are defined via 
the decomposition
$$ \s^{D_{n}}(X_n)= \sum_{w \in W(D_n)}~\s_{w}^{D_{n}}(X_n)~u_{w},$$
$(Da)$~~The type $D$ Schubert polynomials of the second kind are defined via 
the decomposition
$$ \D (X_n)= \sum_{w \in W(D_n)}~\D_{w}(X_n)~u_{w}.$$
\end{dfn}
\begin{rem}{ If we replace in the above formulas the Schubert expressions 
on the Grothendieck expressions, and decomposition in the Nil--Coxeter 
algebras on that in the corresponding $Id$-Coxeter algebras, we will come to 
the  definition of {\it Grothendieck} polynomials of type $B$
 (resp. of types $C$ and $D$) of the first or second kind. 
}
\end{rem}
\begin{lem} (\cite{FK3})~~~~(Factorization formula, the case of type $B$)
$$ \s^{B_{n}}(X_n)=\s^{A_{n-1}}( X_n^{op})~\prod_{i=0}^{n-1} \bigg(
~ h_{0}(x_{n-i})~
 \prod_{j=1}^{n-i-1}~~h_{j}(x_{n-i-j}+x_{n-i}) \biggr) .$$
\end{lem}
\begin{prop} (Factorization formula, the case of type $D$)

$(1)$~~Assume that $ n=2k+1 \ge 3$ is odd, then ~~~~
$\s^{B_n}(X_n)=\s^{A_{n-1}}(X_{n-1}^{op},x_n)~$
$$\prod_{r=k}^{1}~
h_{{\hat 1}}(x_{2r-1}+x_{2r}) \biggl( \prod_{a=2}^{2r}~h_{a}(x_{2r-a}+x_{2r}) 
~\prod_{a=1}^{2r-1}~ h_{a}(x_{2r-1-a}+x_{2r-1})) \biggr).$$
$(2)$~~Assume that $n=2k \ge 4$ is even, then~~~~
$\s^{B_n}(X_n)=\s^{A_{n-1}}(X_{n-1}^{op},x_n)~$
$$ \prod_{r=k-1}^{1}~
h_{{\hat 1}}(x_{2r}+x_{2r+1}) \biggl( \prod_{a=2}^{2r+1}
~h_{a}(x_{2r+1-a}+x_{2r}) 
~\prod_{a=1}^{2r} h_{a}(x_{2r-a}+x_{2r})) \biggr)~h_{{\hat 1}}(x_1),$$
where we set $x_{0}=0,$ ~and $X_n^{op}:=(x_n,\ldots,x_1).$
\end{prop}
Note that the number of terms in each expression is equal to $n(n-1),$ the 
length of the maximal element in the group $W(D_n).$ In fact, the both 
products correspond the ``maximal'' reduced decomposition of the element of 
maximal length in the Weyl group $W(D_n),$ as well as for the $B_n$-case 
stated in Lemma~3.1. 

\subsection{Double Schubert and Grothendieck polynomials}
\begin{dfn} ${}$

$(A)$~~The double Schubert expression $\s^{A_{n-1}}(X,Y)$ and double Schubert 
polynomials $\s^{A_{n-1}}_{w}(X,Y),$~$w \in W(A_{n-1}),$ ~of type
$A$ are defined as follows
$$ \s^{A_{n-1}}(X,Y)=(\s^{A_{n-1}}(-Y))^{-1}~\s^{A_{n-1}}(X)=
\sum_{w \in {\mathbb S}_{n}} \s_{w}(X,Y)~u_{w},$$
$(B)$~~The type $B$ double Schubert expression $\s^{B_{n}}(X,Y)$ and type $B$ 
 double Schubert polynomials of the first kind 
$\s^{B_{n}}_{w}(X,Y),$~$w \in W(B_n),$ are defined as follows
$$ \s^{B_{n}}(X,Y)=(\s^{B_{n}}(-Y))^{-1}~\s^{B_{n}}(X)=
 \sum_{w \in W(B_n)}~\s^{B_{n}}_{w}(X,Y),$$
$(Ba)$~~The type $B$ double Schubert expression $\b (X,Y)$  and type $B$ 
double Schubert polynomials of the second kind 
$\b_{w}(X_n,Y_n),$ ~$w \in W(B_n),$ ~are defined as follows
$$\b(X_n,Y_n)=(\b(-Y_n))^{-1}~\b(X_n)= \sum_{w \in W(B_n)}~\b_{w}(X_n,Y_n).$$
$(C)$~~The type $C$ double Schubert expression $\s^{C_{n}}(X,Y)$ and type $C$ 
 double Schubert polynomials of the first kind 
$\s^{C_{n}}_{w}(X,Y),$~$w \in W(C_n),$ are defined as follows
$$ \s^{C_{n}}(X,Y)=(\s^{C_{n}}(-Y))^{-1}~\s^{C_{n}}(X)=
 \sum_{w \in W(C_n)}~\s^{C_{n}}_{w}(X,Y),$$
$(Ca)$~~The type $C$ double Schubert expression $\c(X_n,Y_n)$  and type $C$ 
double Schubert polynomials of the second kind 
$\c_{w}(X_n,Y_n),$ ~$w \in W(C_n),$ ~are defined as follows
$$\c(X_n,Y_n)=(\c(-Y_n))^{-1}~\c(X_n)= \sum_{w \in W(C_n)}~\c_{w}(X_n,Y_n).$$
$(D)$~~The type $D$ double Schubert expression $\s^{D_{n}}(X,Y)$ and the 
type $D$ double Schubert polynomials $\s^{D_{n}}_{w}(X,Y),$ $w \in W(D_n),$ 
of the first kind are defined as follows
$$\s^{D_{n}}(X,Y)=(\s^{D_n}(-Y))^{-1}~\s^{D_n}(X)=
\sum_{w \in D_n}~ \s^{D_{n}}_{w}(X,Y),$$
$(Da)$~~The type $D$ double Schubert expression $\D(X_{n},Y_{n})$ and the 
type $D$ double Schubert polynomials $\D_{w}(X_n,Y_n),$ $w \in W(D_n),$ 
of the second kind are defined as follows
$$\D(X_{n},Y_{n})=({\D}(-Y_n))^{-1}~{\D}(X_n)=
\sum_{w \in D_n}~ {\D}_{w}(X_n,Y_n).$$
\end{dfn}${}$
\begin{dfn}~~The double Grothendieck polynomials of the first and second types
 are defined by replacing in the above formulas the corresponding double 
Schubert expressions on the corresponding double Grothendieck ones.
\end{dfn}
\begin{lem} (\cite{FS},\cite{FK2})~~The double Schubert expression of the 
type $A$ has the following decomposition
$$\s^{A_{n-1}}(X,Y)= \prod_{i=1}^{n-1}~\prod_{j=n-i}^{1} h_{i+j-1}(x_i+y_j).$$
\end{lem}
\begin{ex} ( $B_2$ double Schubert polynomials of the first kind ) \\
$\s_{id}^{B_2}(X,Y)=1,$  \\
$\s_{u_{0}}^{B_2}(X,Y)=x_1+x_2+y_1+y_2,$ \\
$\s_{u_{1}}^{B_2}(X,Y)=x_1+2x_2+y_1+2y_2,$ \\
$\s_{u_{01}}^{B_2}(X,Y)=x_1~x_2+x_2^2+(y_1+y_2)(x_1+2x_2+y+1+y_2),$ \\
$\s_{u_{10}}^{B_2}(X,Y)=(x_1+x_2)(x_1+x_2+y_1+2 y_2)+y_1~y_2+y_2^2,$ \\
$\s_{u_{010}}^{B_2}(X,Y)= x_1~x_2(x_1+x_2)+(x_1+x_2)(y_1+y_2)(x_1+x_2+y_1+y_2)+y_1~y_2(y_1+y_2),$   \\
$\s_{u_{101}}^{B_2}(X,Y)=(x_1+x_2)~x_2^2+x_2(x_1+x_2)(y_1+2 y_2)+
(x_1+2 x_2)y_2(y_1+y_2)+(y_1+y_2) y_2^2,$ \\
$\s_{u_{0101}}^{B_2}(X,Y)= (x_2+y_2)~\s^{B_{2}}_{u_{010}}(X,Y).$ 
\end{ex}
Let us remark that $\s_{u_{010}}^{B_2}(X,Y)=s_{(2,1)}(X_2,Y_2),$ ~~and~~ 
$\s_{u_{0101}}^{B_2}(X,Y)= P_{(2,1)}(X_2,Y_2).$
\begin{thm}${}$

$(1)$~~(Factorization formula)~~~~$\s^{B_{n}}(X,Y)=$
$$\prod_{i=1}^{n} \biggl(
 \prod_{j=1}^{i-1}~h_{j}(y_{i-j}+y_{i})~~h_{0}(y_{i}) \biggr)~
\s^{A}(X^{op},Y^{op})~\biggl( \prod_{i=n}^{1}~h_{0}(x_{i})~\prod_{j=1}^{i-1}~
h_{j}(x_{i-j}+x_{i}) \biggr).$$
$(2)$~~~~(Specialization formula for $B_n$ double Schubert expression of the 
second kind) \\
${\b}(X_n,- X_n)=(\s^{A_{n-1}}(X_n))^{-1}~B(X_n)~
\s^{A_{n-1}}(-X_n)=$
$$ \prod_{j=1}^{n}~~\biggl(
\prod_{a=j-1}^{1}~h_{a}(x_{a}+x_{j})~h_{0}(x_{j})~\prod_{a=1}^{j-1}~h_{a}(x_{j}-x_{a}) \biggr).$$
\end{thm}
Note that the number of terms in the last product is equal to $n^2,$ i.e. to 
the length of the maximal element in $B_n.$~Recall that by definition
$X_n^{op}=(x_n, \dots, x_1).$ 
  
\subsection{Main properties of $B_n$ double Schubert and Grothendieck 
polynomials}
\begin{thm}~~(The case of double Schubert polynomials)  ${}$

$(1)$~~(Action of divided difference operators)
$$ \partial_{i}^{(x)} \s^{B_{n}}(X,Y)=\s^{B_{n}}(X,Y)~u_{i},
~~~~\partial_{i}^{(y)} \s^{B_{n}}(X,Y)=u_{i}~\s^{B_{n}}(X,Y) $$
for all $i=1,\ldots, n-1,$
$$\partial_{i}^{(x)} \b(X_n,Y_n)=\b(X_n,Y_n)~u_{i},
~~~~\partial_{i}^{(y)} \b(X,Y)=u_{i}~\b(X,Y) $$
for all $i=0,\ldots, n-1;$

$(2)$~~(Stability) ~~~Let $ \iota: B_n \rightarrow B_{n+1}$ be the standard 
embedding of the group $W(B_{n})$ to that $W(B_{n+1}).$ If $w \in W(B_{n}),$ 
then
$$ \s^{B_{n+1}}_{\iota(w)}(X_{n+1},Y_{n+1})=\s^{B_{n}}_{w}(X_{n},Y_{n}),~~~~~
\b_{\iota(w)}(X_{n+1},Y_{n+1})=\b_{w}(X_{n},Y_{n}); $$

$(3)$~~(Vanishing property of double Schubert polynomials of the second 
kind corresponding to a Weyl group $W$ of classical type) \\
Let $w \in W$ and 
$(a_1,\ldots,a_l) \in R(w)$ be a reduced decomposition of the element $w.$ Then
\begin{equation}
\s_{W}(X_{n},-w(X_{n})) =\prod_{r=l}^{1}~h_{a_r}
(x_{s_{a_1} \cdots s_{a_{r-1}}(a_r)} - x_{s_{a_1} \cdots s_{a_r}(a_r)}).
\end{equation}
Therefore, if $v,w \in W,$ then
$$ \b_{w}(X_{n},-v(X_{n})) \not= 0,$$
if and only if $v \leq w$ with respect to the Bruhat order on the group 
 $W.$

$(3a)$~~ Let $w \in {\mathbb S}_n \subset W(B_n),$ then~~~ ${\b}(X_n,-w(X_n))=$
$$ \s^{B_n}(X_n,-w(X_n))=\prod_{(i,j) \in D(w)}~h_{n(i,j)}~(x_{w(i)}-x_j)
=\s^{A_{n-1}}(-X,w(X)),$$
where $D(w)$ denotes the diagram of a  permutation $w \in {\mathbb S}_n,$ see 
e.g. \cite{Mac}, p. 8, and $n(i,j)$ is the number in the box $(i,j) \in D(w)$ 
according to the  standard numbering of the boxes of the diagram $D(w),$ and 
the product is taken according to the reading of the boxes of the diagram 
$D(w)$  column by column, from the bottom to the top, starting from the first 
column, next is the second and so on.

$(3b)$~~Let $w={\bar u} \in W(B_n)$ be the sign permutation corresponding to a 
permutation $u \in {\mathbb S}_n.$ ~~~Then ${\b}(X_n,-w(X_n))=$
$$(\s^{A_{n-1}}(u(X_n)))^{-1}~B(X_n)~\s^{A_{n-1}}(-X_{n})=\prod_{(i,j) \in 
 {\overline D(u)}}~h_{n(i,j)}(x_{w(i)}-x_{j}),$$
where ${{\overline D(u)}}= [1,n]^2 \setminus D(u);$ the product is taken 
according to the reading of the boxes of the set ${\overline D(w)}$  column by 
column, from the bottom to the top, starting from the first column, next is 
the second and so on; by definition we set $x_{{\bar i}}=-x_i.$

$(4)$~~If $w \in W(B_{n})$ ~is a type $B$ Grassmannian permutation of shape 
$\lambda:=\lambda(w),$ then
$$ {\b}_{w}(X_n,Y_n)= P_{\lambda}(X_{l(\lambda)} ~|~ Y_{l(\lambda)}),$$
where $P_{\lambda}(X,Y)$ denotes the factorial Schur polynomial corresponding
 to a partition $\lambda$ introduced and studied in \cite{Iv1}, \cite{Iv2}. 
\end{thm}
{\bf Proof} ${}$

$(1)$~~~By definition, ~~$\partial_{0}^{(x)} {\b}(X_n,Y_n):=x_{1}^{-1}~
 (\b(X_n,Y_n)-{\b}(-x_{1},X_n \setminus \{x_{1} \},Y_n))$=
$$ x_1^{-1}~\biggl( \s^{A_{n-1}}(-Y)^{-1} \sqrt{B(Y)B(X)}~ \biggl( A_{1}(-x_1)
-B(-x_1)~A_{1}(x_1) \biggr)~A_{2}(x_2) \cdots A_{n-1}(x_{n-1}) \biggr) .$$ 
By definition, $B(-x_1)~A_{1}(x_1)= A_{1}(-x_1)~h_{0}(-x_{1}),$ and therefore, 
$A_{1}(-x_1)-B(-x_1)~A_{1}(x_1)= x_1~A_{1}(-x_{1})~u_{0}.$ ~~~Thus, 
$\partial_{0}^{(x)}{\b}(X_n,Y_n)=$ 
$$ x_{1}^{-1} \biggl(\s^{A_{n-1}}(-Y_n) \sqrt{B(Y)B(X)}x_1~A_{1}(-x_{1})~u_{0}
 \biggr)~A_{2}(x_{2}) \cdots A_{n-1}(x_{n-1})= {\b}(X_n,Y_n)~u_{0}.$$ 
Similar reasoning shows that $\partial_{0}^{(y)}~{\b}(X_n,Y_n)= 
u_{0}~{\b}(X_n,Y_n).$ 

{\rm Because the both functions $B(X)~B(Y)$ and $\sqrt{B(X)~B(Y)}$ are 
symmetric with  respect to variables $X$ (resp. $Y$), the divided difference 
operators $\partial_{i}^{(x)}$ (resp. $\partial_{i}^{(y)}$),
~$1 \le i \le n-1,$~ act in fact only on the component $\s^{A_{n-1}}(X)$ 
(resp. $(\s^{A_{n-1}}(-Y))^{-1}$) of the $B_n$ double Schubert expressions 
either the first or the second types. It is well-known  \cite{FS}, \cite{FK2},
 that 
$\partial_{i}^{(x)}~\s^{A_{n-1}}(X)= \s^{A_{n-1}}(X)~u_i$ (resp. 
$\partial_{i}^{(y)} \s^{A_{n-1}}(-Y)^{-1}=u_{i}~ \s^{A_{n-1}}(-Y)^{-1}$).  }

$(2)$~~It is clear that $B_{n+1}(x)=h_n(x)~B_n(x)~h_n(-x),$ and 
$A_{\bullet}^{(n+1)}(x)=h_{n}(x)~A_{\bullet}^{(n)}(x).$  ~Thus if 
$w \in B_n \hookrightarrow B_{n+1},$ then one can erase all appearances of 
the factor $h_n(x)$ in a Schubert expression for $W(B_{n+1})$ to obtain that 
for the group $W(B_n).$ \\
$\qed$  
\begin{rem}~~~{\rm It follows from arguments used in the proof of Theorem~2 
from  \cite{B}, that for any crystallographic Coxeter group $W$  the 
polynomial defined in the Nil--Coxeter algebra $Nil(W)$ by means of the RHS 
of $(8)$ is well-defined, i.e. does not depend on a choice of a reduced 
decomposition of $w \in W$ taken. In fact, one can show that this statement is 
valid for any Coxeter group. One of the main results of our paper states that 
for a Lie group $G$ of classical type, the polynomial in question is equal to 
the specialization $Y=-w(X)$ of the corresponding double Schubert expression 
of the second type, see Definition~3.4. We expect that a similar statement 
can be  generalized for any crystallographic Coxeter group. Let us stress that the equality $(8)$ is true as that among {\it polynomials}, but not only in 
the corresponding equivariant cohomology ring.  \\
}
\end{rem}
\begin{thm}~~(The case of double Grothendieck polynomials)  ${}$

$(1)$~~(Action of isobaric divided difference operators)

$(I)$~~~$\pi_i ^{A}(\g^{A_{n-1}}(X_n))= \g^{A_{n-1}}(X_n)~(u_i - \beta).$
$$ \pi_{i;x}^{(B)}( \g^{B_{n}}(X,Y))=
\g^{B_{n}}(X,Y)~(u_{i}-{\beta \over 2}),
~~~~\pi_{i;y}^{(B)} (\g^{B_{n}}(X,Y))=
(u_{i}-{\beta \over 2})~\g^{B_{n}}(X,Y) $$
for all $i=1,\ldots, n-1;$
$$\pi_{i;x}^{(B)}(\g_{B_n}(X_n,Y_n))=
\g_{B_n}(X_n,Y_n)~(u_{i}-{\beta \over 2}),
~~~~\pi_{i;y}^{(B)}(\g_{B_n}(X,Y))=(u_{i}-{\beta \over 2})~\g_{B_n}(X,Y) $$
for $i=,\ldots, n-1;$

$(II)$ 
$$\pi_{0;x}^{(B)}(\g_{B_n}(X_n,Y_n))=
\g_{B_n}(X_n,Y_n)~(u_{0}-\beta ),
~~~~\pi_{0;y}^{(B)}(\g_{B_n}(X,Y))=(u_{0}-\beta)~\g_{B_n}(X,Y). $$
$$\pi_{{\hat 1};x}^{(D)}(\g_{D_n}(X_n,Y_n))=\g_{D_n}(X_n,Y_n)
~(u_{{\hat 1}}-\beta ),
~~~~\pi_{0;y}^{(D)}(\g_{D_n}(X,Y))=(u_{{\hat 1}}-\beta)~\g_{D_n}(X,Y). $$

$(2)$~~(Stability) ~~~Let $ \iota: B_n \rightarrow B_{n+1}$ be the standard 
embedding of the group $W(B_{n})$ to that $W(B_{n+1}).$ If $w \in W(B_{n}),$ 
then
$$ \s^{B_{n+1}}_{\iota(w)}(X_{n+1},Y_{n+1})=\s^{B_{n}}_{w}(X_{n},Y_{n}),~~~~~
\b_{\iota(w)}(X_{n+1},Y_{n+1})=\b_{w}(X_{n},Y_{n}); $$

$(3)$~~(Vanishing property of double Grothendieck polynomials of the second 
kind corresponding to a Weyl group $W$ of classical type) \\
Let $w \in W$ and 
$(a_1,\ldots,a_l) \in R(w)$ be a reduced decomposition of the element $w.$ Then
\begin{equation}
\g_{W}(X_{n},-w(X_{n})) =\prod_{r=l}^{1}~h_{a_r}
(x_{s_{a_1} \cdots s_{a_{r-1}}(a_r)} -_{\beta} 
x_{s_{a_1} \cdots s_{a_r}(a_r)}).
\end{equation}
Therefore, if $v,w \in W,$ then
$$ \g_{w}(X_{n},-v(X_{n})) \not= 0,$$
if and only if $v \leq w$ with respect to the Bruhat order on the group $W.$
\end{thm}
Before to start the proof, let's state a simple, but useful identity
$$h_{*}(\phi_{\beta}(x))~h_{*}(\phi_{\beta}(-x))=1,$$
where $*$ can be equal to $i, 0 \le i \le n-1,$ or ${\hat 1}.$ 

{\bf Proof}

By definition, $\pi_{i}^{A}(\g^{A_{n-1}}(X_n))=$ \\
 $\pi_{i}^{A}((A_1(x_1) \cdots 
A_{i}(x_i)~A_{i}(x_{i+1})~h_{i}(x_{i+1})^{-1}~A_{i+2}(x_{i+2}) \dots 
A_{n-1}(x_{n-1})) =$ \\
$A_1(x_1) \cdots A_{i}(x_i)~A_{i}(x_{i+1})~\pi_{i}^{A}(h_{i}(x_{i+1})^{-1})~
A_{i+2}(x_{i+2}) \dots A_{n-1}(x_{n-1}).$ Now one can compute~~~
$\pi_{i}^{A}(h_{i}(x_{i+1})^{-1})=$
$$ {1 \over x_i -x_{i+1}} \biggl(
(1+ \beta x_{i+1})(1-{x_{i+1} \over 1+\beta x_{i+1}}~u_{i})-
(1+ \beta x_{i})(1-{x_{i} \over 1+\beta x_{i}}~u_{i}) \biggr)=u_{i}-\beta.$$ 
Therefore, $\pi_i ^{A}(\g^{A_{n-1}}(X_n))= \g^{A_{n-1}}(X_n)~(u_i - \beta).$ \\
Similarly, one can compute $\pi_{i}^{B}~\biggl(h_{i}(\phi_{\beta}(x))^{-1} 
\biggr)=\pi_{i}^{B} (h_{i}(\phi_{\beta}(-x)))= u_{i}-{\beta \over 2}.$ 
Therefore, $\pi_{i;x}^{(B)}( \g^{B_{n}}(X,Y))=
\g^{B_{n}}(X,Y)~(u_{i}-{\beta \over 2}),$~The other cases listed in 
$(I)$ can be proved in a similar manner.~~  $QED.$

Now let's compute the action of isobaric divided difference operator 
$\pi_{0;x}^{B}.$ ~We have~~~ 
$\pi_{0;x}^{(B)}(\g_{B_n}(X_n,Y_n))=$ \\
$ \pi_{0;x}^{(B)} \biggl( 
(\s^{A_{n-1}}(\phi_{\beta}(-Y_n)))^{-1}~\sqrt {{\cal B}(x_1)
 \cdots {\cal B}(x_n)}~\s^{A_{n-1}}(\phi_{\beta}(-X_n)) \biggr)=$ \\
$(\s^{A_{n-1}}(\phi_{\beta}(-Y_n)))^{-1}~\sqrt {{\cal B}(x_2)
 \cdots {\cal B}(x_n)}$ \\
$\biggl(\pi_{0;x}^{(B)}(\sqrt{{\cal B}(x_1)}~{\cal A}_{1}(-x_1)) 
 \biggr)~{\cal A}_{2}(-x_2) \cdots {\cal A}_{n-1}(-x_{n-1}).$ 

Now one can compute ~~~
$ \pi_{0;x}^{(B)}(\sqrt{{\cal B}(x_1)}~{\cal A}_{1}(-x_1))=$ \\
$\sqrt{{\cal B}(x_1)}~{\cal A}_{1}(-x_1) \biggl((1-{\beta \over 2}x_1)-
(1+{\beta \over 2}){\cal A}_{1}(-x_1)^{-1}~{\cal B}(-x_1)~
{\cal A}_{1}(-x_1)\biggr)~x_1^{-1}=$ \\
$\sqrt{{\cal B}(x_1)}~{\cal A}_{1}(-x_1) \biggl((1-{\beta \over 2}x_1)- 
(1+{\beta \over 2})~h_{0}(\phi_{\beta}(-x_1)) \biggr)~x_1^{-1}=$ \\
$\sqrt{{\cal B}(x_1)}~{\cal A}_{1}(-x_1)~(u_0-\beta).$ ~~~Therefore, \\
$\pi_{0;x}^{(B)}(\g_{B_n}(X_n,Y_n))=
\g_{B_n}(X_n,Y_n)~(u_{0}-\beta ).$ Formula for the action of 
$\pi_{0;y}^{(B)}$ can be proved in a similar fashion.~~   $QED$

Finally, let's compute the action of isobaric divided difference operator 
$\pi_{{\hat 1};x}^{D}.$ ~Like in the cases considered before, using the 
identity ~${\cal D}(x)~{\cal D}(-x)=1$~it is enough to compute
~ $\pi_{{\hat 1};x}^{D} \biggl(\sqrt{{\cal D}(x_1)~{\cal D}(x_2)}~
{\cal A}_1(-x_1)~{\cal A}_2(-x_2)\biggr)=$ \\
$\sqrt{{\cal D}(x_1)~{\cal D}(x_2)}~
{\cal A}_1(-x_1)~{\cal A}_2(-x_2)~\biggl((1-{\beta \over 2}x_1)(1-
{\beta \over 2}x_2)-
(1+{\beta \over 2}x_1)(1+{\beta \over 2}x_2)$ \\
${\cal A}_2(-x_2)^{-1}~{\cal A}_1(-x_1)^{-1}~{\cal D}(-x_2)~{\cal D}(-x_1)
{\cal A}_1(x_2)~{\cal A}_2(x_1)\biggr).$ 

Now let's compute the product 
$${\cal A}_2(-x_2)^{-1}~{\cal A}_1(-x_1)^{-1}~{\cal D}(-x_2)~{\cal D}(-x_1)
{\cal A}_1(x_2)~{\cal A}_2(x_1).$$
To accomplish this, we will use the following formulas 
$${\cal D}(x)={\cal A}_1(x)~h_{{\hat 1}}(\phi_{\beta}(x))~{\cal A}_2(-x)^{-1}=
{\cal A}_2(x)~h_{{\hat 1}}(\phi_{\beta}(x))~{\cal A}_1(-x)^{-1}.$$
Thus, the product in question is equal to \\
${\cal A}_2(-x_2)^{-1}~h_{{\hat 1}}(\phi_{\beta}(-x_1))~{\cal A}_2(x_1)^{-1}~
{\cal A}_2(-x_2)~h_{{\hat 1}}(\phi_{\beta}(-x_1))~{\cal A}_2(x_1)=$ \\
$h_2(\phi_{\beta}(x_2))~{\cal A}_3(-x_2)^{-1}~h_{{\hat 1}}(\phi_{\beta}(-x_1))~{\cal A}_3(-x_2)~h_{2}(\phi_{\beta}(-x_2)+_{\beta} \phi_{\beta}(-x_1)) $\\
${\cal A}_3(x_1)^{-1}~h_{{\hat 1}}(\phi_{\beta}(-x_2))~{\cal A}_3(x_1)~h_{2}(\phi_{\beta}(x_1))=$ \\
$h_2(\phi_{\beta}(x_2))~h_{{\hat 1}}(\phi_{\beta}(-x_1))~h_{2}(\phi_{\beta}(-x_2)+_{\beta} \phi_{\beta}(-x_1))~h_{{\hat 1}}(\phi_{\beta}(-x_2))~h_{2}(\phi_{\beta}(-x_1))=$ \\
$h_{{\hat 1}}(\phi_{\beta}(-x_1)+_{\beta} \phi_{\beta}(-x_2).$ The final 
equality follows from the Yang--Baxter relation
$$h_{{\hat 1}}(x)~h_2(x+_{\beta} y)~h_{{\hat 1}}(y)=
h_2(y)~h_{{\hat 1}}(x+_{\beta} y)~h_2(x).$$
Substituting the value of the polynomial in question to the above calculations,
one can find \\
$\pi_{i;x}^{(B)}(\g_{B_n}(X_n,Y_n))= \g_{B_n}(X_n,Y_n)~ \biggl\{ 
{1 \over x_1+c_2}~ 
\biggl((1-{\beta \over 2}x_1)((1-{\beta \over 2}x_2)-$ \\
$((1+{\beta \over 2}x_1)((1+{\beta \over 2}x_2)~
h_{{\hat 1}}(\phi_{\beta}(-x_1)+_{\beta} \phi_{\beta}(-x_2))\biggr)
 \biggr\}  =$
$\g_{B_n}(X_n,Y_n)~(u_{{\hat 1}}-\beta).$ \\
To deduce the  final equality, 
we have used the following equality 
$$  \phi_{\beta}(-x_1)+_{\beta} \phi_{\beta}(-x_2)=
 -{x_1+x_2 \over (1+{\beta \over 2}x_1)((1+{\beta \over 2}x_2)}.$$
$\qed$
%
\subsection{Main properties of $D_n$-Schubert polynomials}
Let set in this Section $\partial_i:=\partial_i^{D}, ~{\hat 1},1,2,\ldots,n-1.$
\begin{thm}${}$

$(1)$~~(Action of divided difference operators)
$$ \partial_{i}^{(x)} \s^{D_{n}}(X,Y)=\s^{D_{n}}(X,Y)~u_{i},
~~~~\partial_{i}^{(y)} \s^{D_{n}}(X,Y)=u_{i}~\s^{D_{n}}(X,Y) $$
for all $i=1,\ldots, n-1,$
$$\partial_{i}^{(x)} {\D}(X_n,Y_n)={\D}(X_n,Y_n)~u_{i},
~~~~\partial_{i}^{(y)} {\D}(X_n,Y_n)=u_{i}~{\D}(X_n,Y_n) $$
for all $i= {\hat 1},1,2,\ldots, n-1;$

$(2)$~~(Stability) ~~~Let $ \iota: D_n \rightarrow D_{n+1}$ be the standard 
embedding of the group $W(D_{n})$ to that $W(D_{n+1}).$ If $w \in W(D_{n}),$ 
then
$$ \s^{D_{n+1}}_{\iota(w)}(X_{n+1},Y_{n+1})=\s^{D_{n}}_{w}(X_{n},Y_{n}),~~~~~
{\D}_{\iota(w)}(X_{n+1},Y_{n+1})={\D}_{w}(X_{n},Y_{n}).$$
\end{thm}
{\bf Proof} ${}$

By definition ~~~$\partial_{{\bar 1}}^{(x)}~{\D}(X_n,Y_n)=$ \\
${1 \over x_1+x_2}~\bigg( \s^{A_{n-1}}(-Y)^{-1} \sqrt{D(x_n) \cdots 
D(x_2)~D(x_1)}~\s^{A_{n-1}}(-x_1,-x_2,-x_3, \ldots,-x_n)-$ \\
$\s^{A_{n-1}}(-Y)^{-1} \sqrt{D(x_n) \cdots D(-x_2)~D(-x_1)}~
\s^{A_{n-1}}(x_2,x_1,-x_3, \ldots,-x_n) \bigg)=$ \\
$\s^{A_{n-1}}(-Y)^{-1} \sqrt{D(x_n) \cdots 
D(x_2)~D(x_1)}~ A_1(-x_1)~A_2(-x_2) \biggl(1- $ \\
$A_2(-x_2)^{-1}~A_1(-x_1)^{-1}~D(-x_1)~D(-x_2)~A_1(x_2)~A_2(x_1) \biggr)
~A_3(-x_3) \cdots A_{n-1}(-x_{n-1}).$ 

Now let us simplify the expression
$$ A_2(-x_2)^{-1}~A_1(-x_1)^{-1}~D(-x_1)~D(-x_2)~A_1(x_2)~A_2(x_1)$$
using the the following formula from Lemma~2.2 :
$$D(x)~D(y)= A_{1}(x)~A_{1}(y)~h_{{\hat 1}}(x+y)
~h_{1}(-x-y)~A_{1}(-x)^{-1}~A_{1}(-y)^{-1}.$$
 Thus the expression in question is equal to
$$A_2(-x_2)^{-1}~A_{1}(-x_2)~h_{{\hat 1}}(-x_1-x_2)~h_1(x_1+x_2)
~A_1(x_1)^{-1}~A_2(x_1)=h_{{\hat 1}}(-x_1-x_2).$$
Therefore~~~~$\partial_{{\hat 1}}^{(x)}~{\D}(X_n,Y_n)=$ \\
${1 \over x_1+x_2} ~\s^{A_{n-1}}(-Y)^{-1} \sqrt{D(x_n) \cdots 
D(x_2)~D(x_1)} A_1(-x_1)~A_2(-x_2)$ \\
$\bigg( 1-h_{{\hat 1}}(-x_1-x_2) \bigg)~
A_3(-x_3) \cdots A_{n-1}(x_{n-1})= {\D}(X_n,Y_n)~u_{{\hat 1}}.$

\subsection{Schubert polynomials of the third kind}

The type $B$ and type $C$ (double) Schubert expressions of the first kind 
have nice combinatorial properties, but they are not compatible with the 
action of divided difference operators ~$\partial_{0}^{B}$ and 
~$\partial_{0}^{C}$ respectively. 
\begin{dfn} Let us set ~~~~${\widetilde {B_n(x)}}=1+B_n(x),$
${\widetilde {D_n(x)}}=1+D_n(x),$
$${\widetilde {C_n(x)}}=(1+C_n(x))/2=
((A_{1}^{(n)}(x)+A_{1}^{(n)}(-x))/2+x~A_{1}^{(n)}(x)~u_{0})
~A_{1}^{(n)}(-x)^{-1}.$$
\end{dfn}
It is clear that 
$$[{\widetilde {B_n(x)}},{\widetilde {B_n(y)}}]=0,~~
[{\widetilde {C_n(x)}},{\widetilde {C_n(y)}}]=0,
[{\widetilde {D_n(x)}},{\widetilde {D_n(y)}}]=0.$$
\begin{lem} ~~~~$\partial_{0}^{B}~\biggl({\widetilde {B_n(x_1)}}~A_{1}(-x_1) 
\biggr)={\widetilde {B_n(x_1)}}~A_{1}(-x_1)~u_{0},$
$$\partial_{0}^{C}~\biggl({\widetilde {C_n(x_1)}}~A_{1}(-x_1) 
\biggr)={\widetilde {C_n(x_1)}}~A_{1}(-x_1)~u_{0},$$
$$\partial_{{\hat 1}}^{D}~\biggl({\widetilde {D_n(x_1)}}~
{\widetilde {D_n(x_2)}}~A_{1}(-x_1)~A_{2}(-x_2) 
\biggr)={\widetilde {D_n(x_1)}}~{\widetilde {D_n(x_2)}}~A_{1}(-x_1)~
A_{2}(-x_2)~u_{{\hat 1}},$$
\end{lem}
{\bf Proof} ${}$

By definition,~~ $\partial_{0}^{C}~\biggl({\widetilde {C_n(x_1)}}~A_{1}(-x_1)
 \biggr)={1 \over 2x} \biggl( x~A(x)~u_{0}+x~A(-x)~u_{0}\biggr)=
{1 \over 2}~\biggl( A(x)+A(-x)\biggr)~u_{0}={\widetilde {C_n(x_1)}}~A_{1}(-x_1)~u_{0}.$ 
\begin{dfn}~~The $C$ type Schubert expression and Schubert polynomials of the 
third kind \cite{FK3} are given by
$${\mathfrak {c}}(X_{n})={\widetilde {C_n(x_1)}}~{\widetilde {C_n(x_2)}} \cdots
{\widetilde {C_n(x_n)}}~\s^{A_{n-1}}(-X_n)= \sum_{w \in W(C_n)} 
{\mathfrak {c}}_{w}(X_{n})~u_{w}.$$
\end{dfn}
\begin{prop} 
~~~${\mathfrak{c}}_{w}(X_{n}) \in \Z_{\ge 0}[x_{1},\ldots,x_n]$ for 
all $w \in W(C_n).$
\end{prop}
{\bf Proof} ${}$

For each $k \ge 1$ let us introduce  polynomials 
$${\cal B}_k(x_1,\ldots,x_k)= {A_1(x_1)+A_1(-x_1) \over 2}~
\prod_{a=k}^{2}~\prod_{b=1}^{k-1}~h_{b}(x_{a})+
x_1~A_1(x_1)~\prod_{a=k}^{2}~\prod_{b=1}^{k-1}~h_{b}(x_{a})~u_{0},$$
$${\cal C}(x_1, \cdots,x_k)= B(x_k)~\prod_{a=1}^{k-1}~
{\cal B}_k(x_2, \ldots, x_{a+1}).$$
We claim that 
$${\mathfrak {c}}(X_{n})= \prod_{k=2}^{n}~{\cal C}(x_1, \ldots,x_k).$$
Proof is straightforward and leaves to the reader.

\bigskip


\end{document}